\documentclass[10pt]{article}

\newtheorem{theorem}{Theorem}[section]
\newtheorem{definition}{Definition}

\usepackage{epsfig} 

\title{Stabilization of nonlinear systems with semi-quadratic cost}

\author{Sergey Nikitin \\
\thanks{Department of Mathematics and Statistics, Arizona State University, Tempe, AZ 85287-1804  {\tt \small nikitin@asu.edu}}
}
 
\begin{document}

\maketitle

\begin{abstract}

The paper addresses the stabilization of nonlinear systems with semi-quadratic cost: quadratic with respect to controls and nonlinear for state variables. Paper presents the effective new feedback synthesis procedure. The novel feedback design procedure is based on the ideas borrowed from nonlinear optics and the theory of semi-classical asymptotics.
\end{abstract}

\section{Introduction}

Stabilization is one of the central topics of control theory. Moreover, in numerous applications of control theory one needs not only to stabilize a nonlinear system but also to minimize a certain cost function (like, for example, energy cost). These type of problems is well-known in linear system theory and the solution can be found in the class of linear feedbacks \cite{Yakubovich}. It is also closely related to $H_\infty$-control theory \cite{Zames}. The nonlinear generalization of $H_\infty$ theory is presented  in \cite{vanDerSchaft}. This paper in its spirit is similar to \cite{vanDerSchaft}. However, the main topic of this publication is different and devoted to the more narrow subject of the new synthesis of feedback stabilizers that minimize the semi-quadratic cost functional defined as
$$
\int_0^\infty \varepsilon (x) + \langle u, Q(x) u \rangle dt.
$$ 
This paper follows the same ideology as publications \cite{Nikitin},\cite{Nikitin99}, \cite{Ser_book}. That means, we assume that the problem of optimal stabilization is solved locally in some small neighborhood of the target equilibrium. In the majority of real world applications that local solution is obtained in the framework of classical linear system theory \cite{Kwakernaak}. Then we follow the ideas outlined in \cite{Nikitin} and extend the optimal feedback outside the small neighborhood. This process can be interpreted in optical terms: the stabilization problem in the reverse time can be treated as a problem from the nonlinear optic with sources of "light" located on the level set of the "local" Lyapunov function. The propagation of "light" is described by projection into $x$-space the bicharacteristics of the Lagrangian manifold associated with the semi-quadratic functional. This scenario efficiently solves the problem of optimal stabilization. However, on this path we encounter the following pitfalls. The first is related to caustic phenomena. Our "light" propagates in non-uniform "media" (distortion is due to function $\varepsilon(x))$ and we face scattering phenomena that leads to caustics. The second problem is related to the caustics. However, its true nature is not "scattering of light" but the singularity of optimizers related, in particular, to those semi-quadratic functional where the matrix $Q(x)$ is degenerate. Discussion of this problem for linear systems can be found in \cite{Jurdjevic}. To resolve the difficulties related to caustics we employ Maslov canonical operator. However, this approach does not solve completely the problems connected with singular optimizers, and this problem remains out of the scope of this publication.

\vspace{0.1cm}

\vspace{0.1cm}

\section{Preliminaries}
Consider a system
\begin{equation}
\label{system}
\dot x =  f(x) + g(x) u
\end{equation}
where $u$ is the control input; $x$ denotes the state of the system and $x\in {\rm R}^n,$ $ {\rm R}^n$ -- $n$-dimensional linear real space. $f(x),\;\;g(x) $ are vector fields:  
$$
\forall x \in {\rm R}^{n} \;\;f(x)\in {\rm R}^n \;\;\mbox{ and } \;\;\;g(x)\in {\rm R}^n
$$
and $f(x),\;\;g(x)$ are infinitely many times continuously differentiable with respect to $x.$ We write $f,g \in {\rm C^\infty}.$ Throughout the paper we assume that 
$$
f(0)=0
$$
and ${\rm R}^n$ is equipped with the scalar product and $\| x \|$ denotes the magnitude of $x,$ i.e
$\| x \| =\sqrt {\langle x , x \rangle },$ where $\langle x , x \rangle$ is the scalar product of $x$ with itself.

Consider the optimization problem
\begin{equation}
\label{functional}
\int_0^\infty \varepsilon (x) + \langle u, Q(x) u \rangle dt \;\;\to \;\;\inf_u
\end{equation}
subjected to the constraints
\begin{eqnarray}
\label{initialValueProblem}
\dot x (t) &=&  f(x) + g(x) u  ,\nonumber\\
&&\\
x(t_0)&=&x_0, \nonumber
\end{eqnarray}
where $f,g \in {\rm C}^\infty $ and $u(t,x)$ is the feedback. We also assume that the function $\varepsilon(x) \in {\rm C}^\infty ,\;\;\varepsilon(0) =0$ and
$$
\varepsilon(x) > 0 \;\;\;\forall\;x \in {\rm R}^n \setminus \{0\}.
$$
 Moreover, $Q(x)\in {\rm C}^\infty,$
$$
Q(x) = Q^T(x) 
$$
 and   
$$
\langle u, Q(x) u \rangle > 0 \;\;\forall\;\;u\not=0.
$$

\vspace{0.1cm}

The goal of this paper is the synthesis of the optimal feedback low
\begin{equation}
\label{feedback}
u= u(t,x) 
\end{equation}
that solves the minimization problem (\ref{functional}) and at the same time stabilizes the system (\ref{system}) at the origin over ${\rm R}^n.$  Throughout the paper we assume that $u(t,x)$ takes its values from ${\rm R}^m.$
The stabilization is defined as follows.
\begin{definition}
\label{stabilizationDef} 
The system (\ref{system}) is said to be stabilizable by the feedback (\ref{feedback}) if the solutions $x(t,t_0,x_0)$ of the closed loop system (\ref{initialValueProblem}) satisfy the condition
$$
\lim_{t \to \infty } x(t,t_0,x_0) = 0  \;\;\forall\; x_0 \in  {\rm R}^n , \;\;t_0\in {\rm R}  .
$$
and the equilibrium is stable.
\end{definition}

Throughout the paper we assume that the system (\ref{system}) is controllable on ${\rm R}^n$ (see, e.g., \cite{Ser_book}).

\section{Lagrange-Pontryagin's manifold}

Our main assumption is that the semi-quadratic stabilization problem stated in the previous section admits a local solution near the origin. For majority of practical applications this assumption is verified by means of classical linear theory applied to the linearization of the problem (\ref{system}), (\ref{functional}) in the vicinity of the origin. Thus, we have a feedback $w(x)$ and a Lyapunov function $V(x)$ that are "good enough" in a neighborhood of the origin. In other words, there exists a positive number $\delta$ such that 
\begin{equation}
\label{inequality}
\langle \frac{\partial }{\partial x} V(x) , f(x) + w(x) b(x) \rangle + \langle w(x) , Q(x) w(x) \rangle  <0\;\;\forall \;\;x \in \{V(x) = \delta\}
\end{equation}
and the feedback $w(x)$ solves our semi-quadratic stabilization problem in the neighborhood $\{V(x) \le \delta\}.$ 

\vspace{0.1cm}

Our main goal is to find the optimal feedback outside the neighborhood
$$
\{V(x) \le \delta\}.
$$
In order to do that consider the following optimization problem (with free end)
\begin{equation}
\label{func}
B(x_0) = \inf_u \{\int_{t_0}^T \varepsilon (x) + \langle u , Q(x) u \rangle dt;\;\;V(x(T)) = \delta \} 
\end{equation}
subordinated to constraints
\begin{eqnarray}
\label{sys}
\dot x &=& f(x) + g(x) u \nonumber \\
&& \\
x(0) &=& x_0. \nonumber
\end{eqnarray}

The function $B(x_0)$ is known as the Bellman function \cite{Bellman}. If we are able to construct the Bellman function $B(x)$ outside the neighborhood
$$
\{V(x) \le \delta \}
$$ 
then our problem is essentially solved as long as the Bellman function is at least one time differentiable. In order to calculate $B(x)$ we use the Pontryagin's maximum principle \cite{Pontryagin}. 

\vspace{0.1cm}

Consider the Lagrange-Pontryagin's manifold defined by bicharacteristics (this terminology is borrowed from \cite{Maslov}, \cite{Maslov_Fedoriuk}), the solutions of the hamiltonian system:
\begin{eqnarray}
\label{hsys}
\dot x &=& -\frac{\partial }{\partial p } H(x,p), \nonumber\\
&&\\
\dot p &=& \frac{\partial }{\partial x } H(x,p), \nonumber
\end{eqnarray}
where $H(x,p)$ is
\begin{equation}
\label{hamiltonian}
H(x,p) = \varepsilon (x) + \langle p , f(x)  \rangle - \frac14 \langle p , g(x) Q^{-1}(x) g^T(x) p \rangle. 
\end{equation}
All the bicharacteristics are emitted from the Lagrangian manifold associated with the level set
$$
\{V(x) = \delta \}.
$$
In other words, the initial conditions for (\ref{hsys}) are
$$
x(0) = x_0,\;\;p(0) = \lambda(x_0) \cdot \frac{\partial }{\partial x}V(x_0)
$$
where $x_0\in \{V(x) = \delta \}$ and $\lambda(x_0)$ is such that 
$$
H(x_0, \lambda(x_0) \cdot \frac{\partial }{\partial x}V(x_0) ) =0.
$$
The existence of such $\lambda(x_0)$ follows from (\ref{inequality}).
 The union of all such bicharacteristics for all points  $x_0\in \{V(x) = \delta \}$ is the Lagrange-Pontryagin's manifold denoted in the sequel as $\Lambda .$ 

\vspace{0.1cm}

It follows from Pontryagin's maximum principle that any optimal solution for the problem (\ref{func}), (\ref{sys}) can be lifted to the Lagrange-Pontryagin's manifold  $\Lambda.$ Notice, that the singular optimizers are essentially eliminated by the inequality (\ref{inequality}) and the assumption that the inverse $Q^{-1}(x)$ exists for any $x.$ However, "scattering of light" still remains as the cause of the problems related to the caustic phenomena.

\vspace{0.1cm}

Let us introduce the natural projection $P(x,p)=x,$
$$
P:\;\;\Lambda \;\to \; {\rm R}^n.
$$
 The Lagrange-Pontryagin's manifold admits a natural parametrization $(\tau , \xi ),$ where $\tau$ is the parameter along the bicharacteristics and $\xi$ denotes (local) coordinates on the level set $\{V(x) = \delta \}.$ The set of critical values for the natural projection $P$ is called caustic set or simply caustics. We denote the caustic set by $\partial C.$ The Sard's theorem \cite{Golubitsky} implies that $\partial C$  has zero Lebesgue measure in ${\rm R}^n.$

\vspace{0.1cm}

Now consider the function
\begin{equation}
\label{mainFunction}
S(\tau,\xi) = -\int_\gamma \langle p , dx \rangle  
\end{equation}
where $\gamma$ is the bicharacteristic connecting $(0,\;\xi )$ and $(\tau,\;\xi).$

\vspace{0.1cm}

Then the following statement takes place.

\begin{theorem}
\label{Bellman}
The Bellman function $B(x)$ is defined as 
\begin{equation}
\label{regularBellman}
B(x)= \inf\{S(\tau,\xi);\;\;(\tau,\xi) \in P^{-1}(x)\}.
\end{equation}
The derivative $\frac{\partial }{\partial x} B(x) $ exists for any point $x\in {\rm R}^n \setminus \partial C$ where the infimum is realized by a point $(\hat \tau,\hat \xi) \in P^{-1}(x)$ and $S(\hat \tau,\hat \xi) < S(\tau,\xi)$ for all $(\tau,\xi) \in  P^{-1}(x) \setminus (\hat \tau,\hat \xi).$ Moreover,
$$
H(x,\frac{\partial }{\partial x} B(x))=0.
$$
\end{theorem}
{\bf Proof.} If $(\tau, \xi )\in P^{-1}(x_0)$ then
$$
S(\tau,\xi) =  -\int_\gamma \langle p , dx \rangle = -\int_0^T \langle p , \frac{\partial }{\partial p} H(x,p)  \rangle  dt =
$$
$$
 -\int_0^T H(x,p) -\varepsilon(x) + \frac14 \langle p , g(x) Q^{-1}(x) g^T(x) p \rangle dt
$$
and Pontryagin's maximum principle yields $H(x,p)=0.$
Hence,
$$
S(\tau,\xi) = \int_0^T \varepsilon(x) - \frac14 \langle p , g(x) Q^{-1}(x) g^T(x) p \rangle dt 
$$
and applying once more Pontryagin's maximum principle yields
$$
\inf\{S(\tau,\xi);\;\;(\tau,\xi) \in P^{-1}(x)\} = \inf_u\{ \int_0^T \varepsilon (x) + \langle u , Q(x) u \rangle   dt; \;\; V(x(T))=\delta  \}. 
$$
That proves (\ref{regularBellman}).

\vspace{0.1cm}

Let $x_0 \in {\rm R}^n \setminus \partial C .$ If the infimum is realized by a point $(\hat \tau,\hat \xi) \in P^{-1}(x_0)$ and $S(\hat \tau,\hat \xi) < S(\tau,\xi)$ for all $(\tau,\xi) \in  P^{-1}(x) \setminus (\hat \tau,\hat \xi),$  then there exists a neighborhood $D$ of $ (\hat \tau,\hat \xi)$ such that the natural projection $P$ is a diffeomorphism between $D$ and $P(D).$ Hence,
$$
B(x) = S(\tau(x),\xi(x))
$$
where $(\tau(x),\xi(x))$ denotes the smooth parametrization of $D$ by $x$ from $P(D).$ That yields the existence of $\frac{\partial }{\partial x} B(x) .$ Finally, $H(x,p)=0$ on $\Lambda,$ and consequently
$$
H(x,\frac{\partial }{\partial x} B(x))=0
$$
holds.

\vspace{0.1cm}

{\bf Q.E.D.}

\vspace{0.1cm}

We call optimal control problem (\ref{functional}), (\ref{system}) {\it Bellman-regular} if for any $x$ the infimum (\ref{regularBellman}) is realized by a point $(\hat \tau,\hat \xi) \in P^{-1}(x)$ and $S(\hat \tau,\hat \xi) < S(\tau,\xi)$ for all $(\tau,\xi) \in  P^{-1}(x) \setminus (\hat \tau,\hat \xi).$

\vspace{0.1cm}

 For a Bellman-regular optimal control problem the derivative $\frac{\partial }{\partial x} B(x) $ exists for any point  $x$ where  the infimum (\ref{regularBellman}) is given by a non critical point of the natural projection $P.$ Due to Sard's theorem \cite{Golubitsky} the Bellman function almost everywhere differentiable on ${\rm R}^n$ and consequently the optimal feedback
$$
u = -\frac12 \cdot  Q^{-1} (x)g^T(x) \frac{\partial }{\partial x} B(x)
$$
is defined almost everywhere on ${\rm R}^n.$ However, for practical applications of control theory the feedback has to be defined at each point from ${\rm R}^n.$ In order to achieve this goal we construct a regularization of the Bellman function near caustics with the ideas from Maslov canonical operator \cite{Maslov}, \cite{Maslov_Fedoriuk}.

Let $Cr(B) \subset \Lambda $ denote the set of critical points for $P$ such that the infimum (\ref{regularBellman}) is realized at each point from  $Cr(B).$ Evidently $P(Cr(B))$ is a subset of caustic $\partial C$ where the Bellman function is not differentiable.

\section{Canonical Maslov Operator}
This section outlines the essential properties of the canonical Maslov operator that will be employed in our optimal feedback design. The reader interested in the details is directed to the publications \cite{Maslov}, \cite{Maslov_Fedoriuk}.

\vspace{0.1cm}

We outline the construction of canonical Maslov operator on a Lagrange manifold $\cal L.$

Consider a partition of set $\{1,\;2,\;\dots , \; n\}$ into two subsets $\alpha$ and $\beta$ such that $\alpha \cap \beta = \emptyset $ and $\alpha \cup \beta = \{1,\;2,\;\dots , \; n\}.$ The linear space defined by coordinates $(x_\alpha , p_\beta )$ is called Lagrangian coordinate space (or Lagrangian coordinates) in ${\rm R}^{2n}_{x,p}.$ It is known \cite{Maslov_Fedoriuk} that any point $\lambda \in {\cal L}$ has a neighborhood that can be parametrized by Lagrangian coordinates $(x_\alpha , p_\beta )$ and the Jacobian
\begin{equation}
\label{Jacobian}
  \frac{\partial(x_\alpha , p_\beta )}{\partial (\tau, \xi )} \not= 0
\end{equation}
in this neighborhood.

\vspace{0.1cm}

There exists \cite{Maslov_Fedoriuk} the countable open cover $\{\Omega_j\}_j$ for Lagrange manifold $\cal L ,$  such that each $\Omega_j$ admits Lagrangian coordinates $(x_\alpha , p_\beta ).$ Let $\{ e_j (\lambda ) \}_j$ be a partition of unity subordinate to the cover $\{\Omega_j\}_j.$ That means each $ e_j (\lambda ) $ is an infinitely many times differentiable real function on  $\cal L $ with compact support, $supp( e_j) \subset \Omega_j$. Moreover, each function takes non-negative real values and
\begin{equation}
\label{theUnit}
\sum_j e_j (\lambda ) =1 \;\;\;\forall \; \lambda \in \cal L.
\end{equation}
For each $\Omega_j$ one can construct the generating function $S_j(x_\alpha , p_\beta )$ where $(x_\alpha , p_\beta )$ are the Lagrangian coordinates on $\Omega_j.$ The generating function $S_j(x_\alpha , p_\beta )$ is defined as
$$
S_j(x_\alpha , p_\beta ) = -\langle x_\beta (x_\alpha , p_\beta ) , p_\beta \rangle +   \int_\gamma \langle p , dx \rangle  
$$
where $\gamma$ denotes the bicharacteristic connecting $(0,\;\xi )$ and $(\tau,\;\xi).$

\vspace{0.1cm}

Notice, that
$$
S_j(x) = \int_\gamma \langle p , dx \rangle  
$$
when the local coordinates of $\Omega_j$ are equal $x,\;\;\alpha = \{1,\;2,\;\dots,\; n\}$ and $\beta = \emptyset .$

\vspace{0.1cm}

Consider a function $\varphi (\lambda )$ such that $supp (\varphi) \subset \Omega_j$ and $\beta \not= \emptyset.$ Then one can define the sub-canonical Maslov operator $Mas_k(\Omega_j)$ as follows
\begin{equation}
\label{subCanonical}
Mas_k(\Omega_j) \varphi = (\frac{k}{-2\pi i})^\frac{\vert \beta \vert}{2} \int_{{\rm R}^{\vert \beta \vert}_{p_\beta } } J_j^{-\frac12} e^{ik (S_j(x_\alpha , p_\beta ) + \langle x_\beta ,  p_\beta \rangle )} \varphi (x_\alpha , p_\beta ) d p_\beta,
\end{equation}
where $J_j$ denotes the absolute value of the Jacobian (\ref{Jacobian}) and $\vert \beta \vert $ is the cardinality of the set $\beta ;$
\begin{equation}
\label{absValueOfJac}
J_j= \vert \frac{\partial(x_\alpha , p_\beta )}{\partial (\tau, \xi )} \vert.
\end{equation}

\vspace{0.1cm}

If $\beta = \emptyset,$ then the sub-canonical Maslov operator $Mas_k(\Omega_j)$ is given by
\begin{equation}
\label{regularC}
Mas_k(\Omega_j) \varphi =  J_j^{-\frac12} e^{ik S_j(x)} \varphi (x),
\end{equation}
where
$$
S_j(x) = \int_\gamma \langle p , dx \rangle  .
$$

\vspace{0.1cm}

The canonical Maslov operator for $\cal L$ is constructed out of sub-canonical Maslov operators (\ref{subCanonical}), (\ref{regularC}). In order to do that, one needs to introduce the index
$$
\nu (\Omega_j,\Omega_k) = \big(idx(\frac{\partial x_\beta (x_\alpha , p_\beta)}{\partial p_\beta }) -  idx(\frac{\partial x_{\hat{\beta}} (x_{\hat{\alpha}} , p_{\hat{\beta}})}{\partial p_{\hat{\beta}} })\big) ({\rm mod \; 4} ),
$$
where $(x_\alpha , p_\beta)$ and $(x_{\hat{\alpha}} , p_{\hat{\beta}})$ are coordinates of the same point $\lambda \in \Omega_j \cap \Omega_k;$ and
$$
idx(\frac{\partial x_\beta (x_\alpha , p_\beta)}{\partial p_\beta })
$$
denotes the number of negative eigenvalues of
$$
\frac{\partial x_\beta (x_\alpha , p_\beta)}{\partial p_\beta }.
$$
If $ \Omega_j \cap \Omega_k =\emptyset ,$ then index $\nu (\Omega_j,\Omega_k)$ is set to be equal to zero,
$$
\nu (\Omega_j,\Omega_k) =0 \;\;\mbox{ for }\;\;\Omega_j \cap \Omega_k =\emptyset .
$$

Given a chain $(\Omega_{j_0},\;\Omega_{j_1},\;\dots ,\; \Omega_{j_k} )$ such that $\Omega_{j_{s}}\cap \Omega_{j_{s+1}} \not= \emptyset$ (for $s$ taking values $0,\;\dots,\; k-1)$ one can define the index
$$
\nu (\Omega_{j_0},\;\Omega_{j_1},\;\dots ,\; \Omega_{j_k} ) = \big( \sum_{s=0}^{k-1} \nu (\Omega_{j_{s}}, \Omega_{j_{s+1}})  \big) ({\rm mod \; 4} ).
$$
A chain $(\Omega_{j_0},\;\Omega_{j_1},\;\dots ,\; \Omega_{j_k} )$ is called closed when $\Omega_{j_0} = \Omega_{j_k} .$ We say that Lagrange manifold $\cal L$ admits quantization if the following two conditions hold.
\begin{itemize}
\item[i.] For any closed curve (loop) $\gamma \subset \cal L$ we have
$$
\int_\gamma \langle p , dx \rangle = 0.
$$
\item[ii.] For any closed chain $(\Omega_{j_0},\;\Omega_{j_1},\;\dots ,\; \Omega_{j_k} )$ we have
$$
\nu (\Omega_{j_0},\;\Omega_{j_1},\;\dots ,\; \Omega_{j_k} ) = 0 .
$$
\end{itemize}
It is shown in \cite{Maslov_Fedoriuk} that any simply-connected Lagrangian manifold admits quantization.

\vspace{0.1cm}

Now we can introduce the canonical Maslov operator for Lagrange manifold $\cal L.$ If $\cal L$ admits quantization then the canonical Maslov operator $Mas_k(\cal L)$ is defined as
\begin{equation}
\label{canonical}
Mas_k({\cal L}) \varphi = \sum_j c_j \cdot Mas_k(\Omega_j) (e_j\cdot \varphi ),\;\;\varphi \in C^\infty ({\cal L}),
\end{equation}
where $e_j(\lambda)$ is from (\ref{theUnit}) and $Mas_k(\Omega_j)$ is a subcanonical Maslov operator introduced by (\ref{subCanonical}), (\ref{regularC}). In order to define constants $\{c_j\}$ in (\ref{canonical}) we fix $\Omega_{j_0}$ with coordinates $x$ and choose
$$
c_j = e^{-\frac{i\pi}{2}\nu_j},
$$
where $\nu_j$ denotes the index of a chain connecting $\Omega_{j_0}$ with $\Omega_j.$ Due to quantization of $\cal L$ the index $\nu_j$ does not depend on the choice of the chain. It follows from \cite{Maslov_Fedoriuk} that
$$
Mas_k({\cal L}) :\;\;C^\infty({\cal L}) \;\to \; C^\infty(P({\cal L})),
$$
where $P$ denotes the natural projection $ P(x,p)=x.$ For the sake of convenience we will use $Mas_k({\cal L}) $ instead of $Mas_k({\cal L}) 1.$

\section{Asymptotics of optimal feedback}

We assume that the control problem (\ref{functional}), (\ref{system}) is Bellman-regular. In order to create a regularization for Bellman function we take a countable cover $\{\Omega_j\}$ of $Cr(B)$ such that every $\Omega_j$ has Lagrangian coordinates $(x_{\alpha_j}, p_{\beta_j}).$ Moreover, if a point $(\hat \tau , \hat \xi) \in Cr(B)$ belongs to $\Omega_j$ then
$$
\Omega_j \cap (P^{-1} (P(\hat \tau , \hat \xi))\setminus (\hat \tau , \hat \xi)) = \emptyset.
$$

 The existence of such cover follows from \cite{Maslov_Fedoriuk}. Let us assume that every connected component $\cal L_\ell$ of
$$
\cup_j \Omega_j = \cup_\ell {\cal L}_\ell
$$
is a Lagrangian manifold that admits quantization. Then we have smooth functions $\{Mas_k( {\cal L}_\ell)\}$ defined on each $\{P({\cal L}_\ell)\},$ respectively. Now consider a countable family of open sets
$$
 b(x_\mu,\varepsilon_\mu) = \{\|x - x_\mu\|<\varepsilon_\mu \} ,\;\;d(x_\mu,\frac{\varepsilon_\mu}{2}) = \{\|x - x_\mu\|>\frac{\varepsilon_\mu}{2} \}
$$
such that
$$
P(Cr(B)) \subset \cup_\mu b(x_\mu,\varepsilon_\mu  ),\;\;P(\Lambda) \subset \cup_\mu (b(x_\mu,\varepsilon_\mu  ) \cup d(x_\mu,\frac{\varepsilon_\mu}{2}))
$$
and each $x_\mu \in P(Cr(B)).$ Moreover, the real numbers $\{ \varepsilon_\mu \}$ are chosen so that for each $\mu$ one can find $j$ such that
$$
b(x_\mu,\varepsilon_\mu  ) \subset P(\Omega_j)
$$
and $\Omega_j $ has Lagrangian coordinates.

 The cover evidently exists (see, e.g. \cite{Golubitsky}). Consider partition of unity $\{w_\mu, v_\mu\}$ subordinate to cover
$$
\{b(x_\mu,\varepsilon_\mu), d(x_\mu,\frac{\varepsilon_\mu}{2}) \},
$$
where $\{w_\mu, v_\mu\}$ are infinitely differentiable non-negative functions with compact supports,
$$
supp(w_\mu) \subset b(x_\mu,\varepsilon_\mu), \;\;\;\; supp(v_\mu) \subset d(x_\mu,\frac{\varepsilon_\mu}{2})
$$
and
$$
\sum_\mu (w_\mu + v_\mu) =1 \;\;\forall\;\;x\in {\rm R}^n.
$$
We propose the following regularization for the Bellman function
$$
B(x,k) = B(x) \cdot (\sum_\mu v_\mu) - Re\{\sum_j   \frac{1}{ik\cdot p_0} \ln(Mas_k({\cal L}_j))\} \cdot ( \sum_{x_\mu \in P(Cr(B) \cap {\cal L}_j)}  w_\mu ),
$$
where $i=\sqrt{-1};\;\;\ln(z)$ is a fixed branch of the complex logarithm; $Re(a+ib)=a;\;\;\;B(x)$ denotes the Bellman function constructed in (\ref{regularBellman}). We show that
$$
B(x,k) \;\to\;B(x)\;\;\mbox{ as }\;\;k\;  \to \; \infty
$$
uniformly on any compact set $K$ in $C^1$-topology. Indeed,
$$
B(x) = B(x,k) \;\;\;\mbox{ for } \;\;\;x \notin \cup_{x_\mu \in P(Cr(B)) } b(x_\mu,\varepsilon_\mu).
$$
A compact $K$ intersects only a finite number of balls $\{ b(x_\mu,\varepsilon_\mu)\}$ and we need to prove uniform convergence on closure of each ball. The latter follows from the theory of semi-classical asymptotics \cite{Maslov_Fedoriuk}. Since
$$
b(x_\mu,\varepsilon_\mu  ) \subset P(\Omega_j)
$$
the corresponding function $Mas_k({\cal L}_j)$ takes form
$$
Mas_k({\cal L}_j) =  (\frac{k}{-2\pi i})^\frac{\vert \beta \vert}{2} \int_{{\rm R}^{\vert \beta \vert}_{p_\beta } } J_j^{-\frac12} e^{ik (S_j(x_\alpha , p_\beta ) + \langle x_\beta ,  p_\beta \rangle )} d p_\beta\;\;\;\forall\;\;x\in b(x_\mu,\varepsilon_\mu  ),
$$
where $J_j$ is defined in (\ref{absValueOfJac}). The method of stationary phase (see, e.g. \cite{Guillemin}) yields
$$
Mas_k({\cal L}_j) = C_j J_j^{-\frac12} \cdot \vert \frac{\partial x_{\beta}}{\partial p_\beta} \vert^{-\frac12} e^{-ik\cdot p_0B(x)} + O(\frac{1}{k}) \;\;\mbox{ as } \;\;k\to \infty
$$
where $C_j$ is a constant. This asymptotic can be differentiated and the convergence is uniform with respect to $x.$ Hence,
\begin{equation}
\label{Masymptotic}
\frac{1}{ik \cdot p_0}\ln (Mas_k({\cal L}_j)) = - B(x) + O(\frac{1}{k}) \;\;\;\mbox{ as } \;\;\;k\to \infty
\end{equation}
Thus, we justified the convergence
$$
B(x,k) \;\;\to \;\; B(x)
$$
in $C^1$ - topology over any compact $K.$

 By construction $B_k(x,k)$ is differentiable at any point $x,$ and therefore, can be used for feedback design instead of $B(x).$ Our "semi-classical" feedback
$$
u = -\frac12 \cdot  Q^{-1} (x)g^T(x) \frac{\partial }{\partial x}  B(x,k)
$$
is equal to the original optimal feedback
$$
u = -\frac12 \cdot  Q^{-1} (x)g^T(x) \frac{\partial }{\partial x}  B(x)
$$
everywhere outside an $\{\varepsilon_\mu \}$-neighborhood of $P(Cr(B))$ and the neighborhood can be made as small as necessary for a particular application.

\vspace{0.1cm}

The feedback low can be designed in the same way as it was done in \cite{Nikitin}, \cite{Nikitin99}, \cite{Ser_book}:

$$
u(x) = \left\{ \begin{array}{cc}
                   w(x) & \mbox{ for } x\in \{V(x) \le \delta \} \\
                   -\frac12 \cdot  Q^{-1} (x)g^T(x) \frac{\partial }{\partial x}  B(x,k)& \mbox{ for } x\in \{V(x) > \delta \} 
               \end{array}
       \right.
$$
where $w(x)$ and $V(x)$ are defined in (\ref{inequality}).

\vspace{0.1cm}

We summarize the discussion of this section in the form of the following theorem.

\begin{theorem} Consider  Bellman-regular  control  problem (\ref{system}), (\ref{functional}) that admits optimal solution for its linearization
$$
\dot x = Ax + Bu,\;\; \int_0^\infty \langle x, \frac{\partial^2}{\partial x^2} \varepsilon (0) x \rangle + <u, Q(0) u> dt \;\;\to\;\;\inf_u
$$
where $\varepsilon (0)=0,$ $\frac{\partial}{\partial x} \varepsilon (0) =0$ and $\frac{\partial^2}{\partial x^2} \varepsilon (0)$ denotes the Hessian of the function $\varepsilon (x)$ (in this context $\langle x, \frac{\partial^2}{\partial x^2} \varepsilon (0) x \rangle $ is strictly positive quadratic form); $A=\frac{\partial}{\partial x} f(0)$ and $B=\frac{\partial}{\partial x}g(0).$ 
Then there exists the Bellman function $B(x)\in C^\infty ({\rm R}^n \setminus P(Cr(B)))$ and $B(x)$ is at least one time continuously differentiable in a neighborhood of any point from $ P(Cr(B)).$
\end{theorem}

\vspace{0.1cm}

{\bf Proof.} Since the control  problem (\ref{system}), (\ref{functional}) admits optimal solution for its linearization we have a linear feedback $w(x)$ and a quadratic Bellman function ${\cal W}(x)$ for the linearization of the problem. Moreover, there exists a positive number $\delta$ such that
\begin{equation}
\label{linequality}
\langle \frac{\partial }{\partial x} {\cal W}(x) , f(x) + w(x) b(x) \rangle + \langle w(x) , Q(x) w(x) \rangle  <0\;\;\forall \;\;x \in \{{\cal W}(x) = \delta\}.
\end{equation}
Making use of Theorem \ref{Bellman} yields the existence of the Bellman function $B_\delta(x)$ outside the neighborhood
$$
 \{{\cal W}(x) \le \delta\}
$$
such that 
 $$
 \frac{\partial }{\partial x} B_\delta(x) = p \;\;\forall \;\;x \in P(\Lambda_\delta) \setminus P(Cr(B_\delta))
$$
where $(x,p)$ belongs to Lagrange-Pontryagin's manifold $\Lambda_\delta, $ where $\delta$ is from (\ref{linequality}).  Lagrange-Pontryagin's manifold is continuously differentiable infinitely many times then so is the Bellman function $B_\delta (x)$ for all $x\in  P(\Lambda_\delta) \setminus P(Cr(B_\delta))$ and $x\not=0.$

\vspace{0.1cm}

If $x\in P(Cr(B_\delta))$ then one can construct $C^\infty$-approximation (\ref{Masymptotic})  for $B_\delta (x)$ with the help of the canonical Maslov operator. This approximation converges uniformly in $C^1$-topology to $B_\delta(x).$ Hence, $B_\delta(x)$ is at least one time continuously differentiable at  $x\in P(Cr(B_\delta)).$

\vspace{0.1cm}

Taking the limit $\delta\;\to\;0$ yields that the Bellman function
$$
B(x) = \lim_{\delta\;\to\;0} B_\delta (x).
$$ 
That justifies the assertion of the theorem and the proof is completed.

\vspace{0.1cm}

{\bf Q.E.D.}

\end{document}